\documentclass[10pt,leqno]{amsart} %%%%leqno es para numeraci\'{o}n de
   \usepackage[centertags]{amsmath}
   \usepackage{amsfonts}
   \usepackage{amsmath}
   \usepackage{amssymb}
   \usepackage{amsthm}
   \usepackage{newlfont}
   \usepackage[latin1]{inputenc}%%%%para que reconozca las tildes.
\usepackage{xcolor}

\allowdisplaybreaks[3]
\theoremstyle{plain}
\newtheorem{Theo}{Theorem}[section]

\newtheorem{Cor}[Theo]{Corollary}
\newtheorem{Lem}[Theo]{Lemma}

\theoremstyle{definition}
\newtheorem{Def}{Definition}[section]

\theoremstyle{remark}
\newtheorem{Rem}{Remark}

\newcommand{\hh}{\mathfrak{h}}

\newcommand{\RR}{\mathbb{R}}
\newcommand{\ZZ}{\mathbb{Z}}
\newcommand{\NN}{\mathbb{N}}

\newcommand{\nr}{\operatorname{nr}}

\numberwithin{equation}{section}

\newcommand{\card}{\operatorname{card}}

\newcommand\Hh{{\mathcal H}}

\newcommand\Be{\mathfrak{b}}

\newcommand\sign{\operatorname{sign}}

\catcode`\,12

\title{Zeros of linear combinations of Hermite polynomials}
\author{Antonio J. Dur\'an}
\address{Departamento de An\'a\-li\-sis Mate\-m\'a\-ti\-co and IMUS,
        Universidad de Sevilla,
        41080 Sevilla, Spain}
\email{duran@us.es}
   \date{}

   \thanks{This research was partially supported by PID2021-124332NB-C21
(Minis\-te\-rio de Cien\-cia e Inno\-va\-ci\'on and Feder Funds (European Union)), and
FQM-262 (Jun\-ta de Anda\-lu\-c\'ia).}



\begin{document}
   \maketitle

\begin{abstract}
We study the number of real zeros of finite combinations of $K+1$ consecutive normalized Hermite polynomials of the form
$$
q_n(x)=\sum_{j=0}^K\gamma_j\tilde H_{n-j}(x),\quad n\ge K,
$$
where $\gamma_j$, $j=0,\dots ,K$, are real numbers with $\gamma_0=1$, $\gamma_K\not =0$. We consider two different normalizations of Hermite polynomials: the standard one (i.e. $\tilde H_n=H_n$),
and $\tilde H_n=H_n/(2^nn!)$ (so that $q_n$ are Appell polynomials: $q_n'=q_{n-1}$). In both cases, we show the key role played by the polynomial $P(x)=\sum_{j=0}^K\gamma_jx^{K-j}$ to solve this problem.
In particular, if all the zeros of $P$ are real then all the zeros of $q_n$, $n\ge K$, are also real.
\end{abstract}

\section{Introduction}
In the recent paper \cite{Dur0}, we have proved that for any positive measure $\mu$ in the real line, having moments of any order and infinitely many points in its support, there always exists a sequence of
orthogonal polynomials $(p_n)_n$ with respect to $\mu$ such that for any positive integer $K$ and any $K+1$ real numbers
$\gamma_j$, $j=0,\dots ,K$, with $\gamma_0=1$, $\gamma_K\not =0$, the polynomial
\begin{equation}\label{dqnv}
q_n(x)=\sum_{j=0}^K\gamma_jp_{n-j}(x),\quad n\ge K,
\end{equation}
has only real zeros for $n$ big enough (depending on $K$ and the $\gamma_j$'s).
Shohat \cite{Sho} was probably the first to observe that the orthogonality of the sequence $(p_n)_n$ implies that $q_n$ has at least $n-K$ real zeros in the convex hull of the support of $\mu$ (using the usual proof that $p_n$ has its $n$ zeros in the convex hull of the support of $\mu$).
Some other related results on zeros of linear combinations of the form (\ref{dqnv}) can be found in \cite{Peh1,Peh2,IsSa,IsNo,INS,BDr,DJM,KRZ}.

We have to notice that the problem of studying the zeros of finite linear combinations of orthogonal polynomials of the form (\ref{dqnv})
is strongly dependent on the normalization of the polynomials $(p_n)_n$.
We have also proved in \cite{Dur0} that our result applies to the usual normalization of the Hermite polynomials.
The purpose of this paper is to show that the spectral properties of Hermite polynomials allow to prove some more interesting results on the zeros of finite linear combinations of two different normalizations of this classical family of orthogonal polynomials: the standard one $H_n$ and the normalized Hermite polynomials  $H_n/(2^nn!)$.

The content of the paper is as follows. In Section \ref{secH2}, we consider the standard normalization of the Hermite polynomials. The starting point is the following result proved in \cite{Dur0}.

\begin{Cor}[Remark 1 in Section 4.1 of \cite{Dur0}]\label{1her}
For any positive integer $K$ and any finite set of $K+1$ real numbers $\gamma_j$, $j=0,\dots, K$, with $\gamma_0=1$ and $\gamma_K\not=0$, the polynomial
\begin{equation}\label{polhn}
q_n(x)=\sum_{j=0}^K\gamma_jH_{n-j}(x)
\end{equation}
has only real zeros for $n\ge \max\{(K-1)^24^{K-2}\max^2\{|\gamma_j|,2\le j\le K\},2K\}$. Moreover the zeros are simple and interlace the zeros of $H_{n-1}$.
\end{Cor}

However, Corollary \ref{1her} can be improved using the spectral property of the Hermite polynomials with respect to its backward shift operator:
\begin{equation}\label{opfi}
\mbox{if $\Lambda f(x)=2xf(x)-f'(x)$, then $\Lambda H_n=H_{n+1}$.}
\end{equation}
Indeed, we improve Corollary \ref{1her} showing that (1) the real-rootedness of the polynomials $q_n$ are strongly dependent on the zeros of the polynomial
\begin{equation}\label{dP}
P(x)=\sum_{j=0}^K\gamma_jx^{K-j},
\end{equation}
and (2) the existence of interlacing properties between the zeros of $q_{n+1}$ and $q_n$ (for the definition of the interlacing property see Definition \ref{dip} below).

\begin{Theo}\label{hqj} If the polynomial $P$ (\ref{dP}) has only real zeros then all the zeros of the polynomial $q_n$ (\ref{polhn}) are real and simple for $n\ge K$, and the zeros of $q_{n+1}$ interlace the zeros of $q_n$. If $P$ has non real zeros, then there exists a positive integer $n_0$, depending only on the non real zeros of $P$ and $K$, such that for $n\ge n_0$ all the zeros of the polynomial $q_n$  are real and simple and the zeros of $q_{n+1}$ interlace the zeros of $q_n$.
\end{Theo}

In order to prove Theorem \ref{hqj}, in Section \ref{secH} using the first order differential operator (\ref{opfi}) and two sequences of real numbers $(\phi_i)_{i\ge 1}$ and $(\psi_i)_{i\ge 1}$,  we introduce a generalization of the Hermite polynomials whose zeros behave nicely.  These polynomials also satisfy other interesting properties such as Turán type inequalities (see Theorem \ref{ttur}). The polynomials $(q_n)_n$ (\ref{polhn}) are particular cases of these generalized Hermite polynomials. We will also prove that when $\psi_i=0$, $i\ge 1$, the class of all generalized Hermite polynomials is the same class as that of all (tipe II) multiple Hermite polynomials (see Remark \ref{mher}).

Iserles, N{\o}rsett and Saff \cite{IsSa,IsNo,INS} were probably the first to point out the key role
of the real zeros of the polynomial $\sum_{j=0}^K \gamma_jx^j$  to prove the real rootedness of the polynomial
$$
\sum_{j=0}^K\gamma_jH_j(x).
$$
In particular, the case of $P$ having only real zeros in Theorem \ref{hqj} recovers
results by Iserles and Saff \cite[Proposition 1]{IsSa}, although they proved it using an approach different to our method.

We have also studied the following normalization of the Hermite polynomials (Section \ref{secH3})
\begin{equation}\label{henv}
\Hh_n(x)=\frac{1}{2^nn!}H_n(x).
\end{equation}
With this normalization the polynomials
\begin{equation}\label{dqnh}
q_n(x)=\sum_{j=0}^K\gamma_j\Hh_{n-j}(x),\quad n\ge K,
\end{equation}
form an Appell sequence, in the sense that they satisfy $q_n'=q_{n-1}$.

In this case,  we prove that again their real-rootedness are strongly dependent on the zeros of the polynomial $P$ (\ref{dP}), although in a different way as it happens when linear combinations of Hermite polynomials of the form (\ref{polhn}) are considered.

\begin{Theo}\label{zeri} Assume that the polynomial $P$ (\ref{dP}) has $N^{\nr}$ non real zeros. Then
\begin{enumerate}
\item The polynomials $q_n$ (\ref{dqnh}) has only real zeros for all $n\ge 0$ if and only if $N^{\nr}=0$. Moreover, the zeros of $q_n$ are simple and the zeros of $q_{n+1}$ interlace the zeros of $q_n$.
\item If $N^{\nr}>0$ then there exists a nonnegative integer $n_0$ (which we take it to be the smallest one) such that the polynomial $q_{n_0}$ has exactly $N^{\nr}$ non real zeros. Moreover $q_n$ has exactly $n-N^{\nr}$ real zeros and $N^{\nr}$ non real zeros if and only if $n\ge n_0$, in which case they are simple and the real zeros of $q_{n+1}$ interlace the real zeros of $q_n$.
\end{enumerate}
\end{Theo}

We have also proved that the polynomials $q_n$ (\ref{dqnh}) satisfy the following Turán type inequality:
$$
q_{n-1}^2(x)-q_n(x)q_{n-2}(x)>0,\quad x\in \RR,
$$
where $n\ge 2$ if all the zeros of $P$ are real, and $n$ has to be taken big enough if $P$ has non real zeros.

Finally, we have studied the asymptotic behaviour of the zeros of the polynomials $q_n$. The case (\ref{polhn}) was studied in \cite[Corollary 4.4]{Dur0} and the case (\ref{dqnh}) in Corollary \ref{cas2}.

\section{Preliminaries}
Along this paper, the interlacing property  is defined as follows.

\begin{Def}\label{dip} Given two finite sets $U$ and $V$ of real numbers ordered by size, we say that $U$ strictly interlaces $V$ if $\min U<\min V$ and between any two consecutive elements of any of the two sets there exists one element of the other.
\end{Def}
Observe that if $U$ interlaces $V$, then either $\card(U) = \card(V)$, and then $\max U < \max V$, or $\card(U) = 1 + \card(V)$, and then $\max U > \max V$. Observe also that the interlacing property is not symmetric, due to the condition $\min U < \min V$.

We will use the following version of Obreshkov theorem (see \cite{Branden}).

\begin{Theo}\label{obre}
Let $p$ and $q$ be real polynomials with $\deg p=1+\deg q$.
Then the zeros of $p$ interlace the zeros of $q$ if and only if all the polynomials in the space $\{\mu p(z)+\lambda q(z): \mu,\lambda\in \RR\}$ has only real and simple zeros.
\end{Theo}

The following elementary Lemmas will be useful (they are Lemmas 2.2 and 3.5 of \cite{Dur0}, respectively).

\begin{Lem}\label{lem2} Define from the numbers $A_j$, $j=0,\dots , K$, $A_0,A_K\not =0$, the polynomial $P_A$ as
$$
P_A(x)=\sum_{j=0}^KA_j x^{K-j}.
$$
If  $\theta $ is a zero of $P_A$, we define the polynomial $P_{B}$ and the numbers $B_j$, $j=0,\dots , K-1$, as
$$
P_{B}(x)=\frac{P_{A}(x)}{x-\theta}=\sum_{j=0}^{K-1}B_j x^{K-1-j}.
$$
Then, on the one hand, we have
\begin{equation}\label{id1l2}
A_j=\begin{cases} B_j-\theta B_{j-1},& j=1,\dots, K-1,\\ B_0,& j=0,\\
-\theta B_{K-1},&j=K.\end{cases}
\end{equation}
And on the other hand
\begin{equation}\label{id2l2}
B_{j}=\sum_{i=0}^{j}\theta^iA_{j-i},\quad 1\le j\le K-1.
\end{equation}
\end{Lem}

Using the notation of Lemma \ref{lem2}, given real numbers $B_j$, $0\le j\le K-1$, and $\theta$, we can produce the real numbers $A_j$, $0\le j\le K$, as in (\ref{id1l2}), so that we have for
$$
P_{A,\theta}=\sum_{j=0}^{K}A_jx^{K-j},\quad P_{B}(x)=\sum_{j=0}^{K-1}B_jx^{K-1-j}
$$
the identity
$$
P_{A,\theta}(x)=(x-\theta)P_{B}
$$
(we have included the real number $\theta$ in the notation $P_{A,\theta}$ to stress the dependence of this polynomial on $\theta$).

If we define
\begin{equation}\label{pidco}
q_n^{B}(x)=\sum_{j=0}^{K-1}B_jp_{n-j},\quad q_n^{A;\theta} (x)=\sum_{j=0}^{K}A_jp_{n-j}
\end{equation}
the identity (\ref{id1l2}) straightforwardly gives
\begin{equation}\label{idco}
q_{n+1}^{A;\theta}(x)=q_{n+1}^{B}(x)-\theta q_{n}^{B}(x),\quad n\ge K.
\end{equation}
We then have.

\begin{Lem}\label{enze}
Assume that all the zeros of the polynomials $q_n^{B}$ are real and simple for $n\ge n_0$. Then the following conditions are equivalent.
\begin{enumerate}
\item The zeros of $q_{n+1}^{B}$ interlace the zeros of $q_n^{B}$ for $n\ge n_0$.
\item For all real number $\theta$ the polynomial $q_n^{A;\theta}$ has only real and simple zeros for $n\ge n_0+1$.
\end{enumerate}
Moreover, in that case the zeros of $q_n^A$ interlace the zeros of $q_n^B$.
\end{Lem}

The following Lemma will be also useful (the proof is similar to the usual proof for the Hurwitz's Theorem (see \cite[p. 178]{ahl}) and it is omitted).

\begin{Lem}\label{mamon}
Let $f_n$, $g_n$, $f$ be analytic functions in a region $\Omega$ of the complex plane. Assume $f$ has $N$ non real zeros in $\Omega$ and that
$$
\lim_nf_n(z)=f(z),\quad  \lim_ng_n(z)=zf(z),
$$
uniformly in compact sets of $\Omega$. Then there exists $n_*\in \NN$ such that for $n\ge n_*$ the function $f_n(z)-\theta g_n(z)$ has at least $N$ non real zeros in $\Omega$ for any real number $\theta$.
\end{Lem}

(We stress that the positive integer $n_*$ guaranteed by the Lemma does not depend on the real number $\theta$).

\section{A generalization of Hermite polynomials}\label{secH}
As we wrote in the Introduction, Corollary \ref{1her} can be improved using the backward shift operator for the Hermite polynomials (\cite[p. 251]{KLS}):
\begin{equation}\label{bsh}
H_{n+1}(x)=2xH_n(x)-H_n'(x).
\end{equation}
We start by using this backward shift operator to construct a generalization of the Hermite polynomials that is interesting in itself.

Let $\Lambda$ be the first order differential operator
\begin{equation}\label{opf}
\Lambda f(x)=2xf(x)-f'(x).
\end{equation}

\begin{Def} A linear operator $T$ acting in the linear space of real polynomials is a real zero increasing operator if for all polynomial $p$ the number of real zeros of $T(p)$ is greater than the number of real zeros of $p$.
\end{Def}

\begin{Lem}\label{hoiz} The operator $\Lambda$ is a real zero increasing operator. Moreover, if all the zeros of $p$ are real and simple then all the zeros of $\Lambda(p)$ are also real and simple, and interlace the zeros of $p$.
\end{Lem}

\begin{proof}
Write $q=\Lambda p=2xp-p'$.

Assume first that $p$ has $k$ real and simple zeros, say
$$
\zeta_1<\dots<\zeta_k.
$$
Then $q(\zeta_i)q(\zeta_{i+1})=p'(\zeta_i)p'(\zeta_{i+1})<0$, and so $q$ has at least one zero in each interval $(\zeta_i,\zeta_{i+1})$, $i=1,\dots, k-1$. It is easy to see that $q$ has also at least one zero in $(-\infty,\zeta_1)$ and $(\zeta_k,+\infty)$. This proves that $q$ has at least $k+1$ zeros. This also shows that if all the zeros of $p$ are real and simple then all the zeros of $q$ are also real and simple, and interlace the zeros of $p$.

If $p$ has zeros of multiplicity bigger than $1$, we can use a continuity argument.

\end{proof}

The property of being the operator $\Lambda$ a real zero increasing operator is, according to the previous Lemma, a naive consequence of its definition. However, Corollary \ref{1her} is actually saying that $\Lambda$ is a real zero increasing operator in the following deeper sense: if $p\not =0$ is any polynomial then there exists $n_0$ (which depends on $p$), such that for all $n\ge n_0$, all the zeros of the polynomial $\Lambda ^n p$ are real, no matter the number of real zeros of $p$. Indeed, given $p\not=0$ of degree $K$, there are real numbers $\gamma_j$ such that
$$
p(x)=\sum_{j=0}^K \gamma_jH_{K-j}(x).
$$
The identity (\ref{bsh}) gives
$$
\Lambda^n p(x)=\sum_{j=0}^K \gamma_jH_{n+K-j}(x).
$$
And hence Corollary \ref{1her} says that for $n$  big enough $\Lambda ^n p$ has only real zeros.

\bigskip

Given two sequences of real numbers $\phi=(\phi_i)_{i\ge1}$ and $\psi=(\psi_i)_{i\ge1}$, we define the sequence of generalized Hermite polynomials associated to $\phi, \psi$ as $\hh_0^{\phi,\psi}=1$ and for $n\ge 1$
\begin{equation}\label{dhg}
\hh_{n+1}^{\phi,\psi}(x)=\Lambda \hh_{n}^{\phi,\psi}(x) +(\phi_{n+1}+x\psi_{n+1})\hh_{n}^{\phi,\psi}(x).
\end{equation}

When $\psi_i\not =-2$, for all $i\ge 1$,  then $\hh_{n}^{\phi,\psi}$ has degree $n$ and leading coefficient equal to $\prod_{i=1}^n(\psi_i+2)$.

In order to simplify the notation, for a real number $u\in \RR$, the sequence $\phi_i=u$, $i\ge 1$, is denoted just by $u$.

As an easy consequence of the backward shift operator for the Hermite polynomials (\ref{bsh}) we have
\begin{equation}\label{rbel}
\hh_n^{0,0}(x)=H_n(x),\quad n\ge 0.
\end{equation}
Moreover, for $r,s\in \RR$, $s\not =-2$, a simple computation gives
$$
\hh_n^{r,s}(x)=\left(\frac{s+2}{2}\right)^{n/2}H_n\left(\left(\frac{s+2}{2}\right)^{1/2}\left(x+\frac{r}{s+2}\right)\right).
$$
In \cite{durgb}, we consider a similar idea to generalize the Bell polynomials
$$
\Be_n(x)=\sum_{j=0}^nS(n,j)x^j,\quad n\ge 0,
$$
(where $S(n,j)$, $0\le j\le n$, denote the Stirling numbers of the second kind) from its backward shift operator
$$
\Be_{n+1}(x)=x\left(1+\frac{d}{dx}\right)\Be_n(x).
$$
This approach has allowed us to show an unexpected connection between Bell and Laguerre polynomials (which we use in the forthcoming paper \cite{Durl} to study the zeros of linear combinations of monic Laguerre polynomials).

\medskip

We next prove that under mild assumption, the zeros of $\hh_n^{\phi,\psi}$ behave nicely.

\begin{Theo}\label{pep} Let $\phi$ and $\psi$ be two sequences of real numbers with $\psi_i>-2$, $i\ge 1$. Then for $n\ge 0$, the polynomial $\hh_n^{\phi,\psi}$ has only real and simple zeros. Moreover, the zeros of $\hh_{n+1}^{\phi,\psi}$ interlace the zeros of $\hh_n^{\phi,\psi}$.
\end{Theo}

\begin{proof}
We proceed by induction on $n$. For $n=0,1$, the result is trivial.

Assume $\hh_n^{\phi,\psi}$ has only real and simple zeros. Write then
$$
\zeta_1<\dots<\zeta_n,
$$
for the zeros of $\hh_n^{\phi,\psi}$.

The definitions (\ref{opf}) and (\ref{dhg}) give
$$
\hh_{n+1}^{\phi,\psi}(x)=2x \hh_{n}^{\phi,\psi}(x)-(\hh_{n}^{\phi,\psi})'(x) +(\phi_{n+1}+x\psi_{n+1})\hh_{n}^{\phi,\psi}(x).
$$
And so
$$
\hh_{n+1}^{\phi,\psi}(\zeta_i)=-(\hh_{n}^{\phi,\psi})'(\zeta_i).
$$
This shows that $\hh_{n+1}^{\phi,\psi}$ has at least one zero in each interval $(\zeta_i,\zeta_{i+1})$, $i=1,\dots, n-1$.

Since the leading coefficient of $\hh_{n+1}^{\phi,\psi}$ is $\prod_{i=1}^{n+1}(\psi_i+2)$ and $\psi_i>-2$, we deduce that the leading coefficients of $\hh_{n+1}^{\phi,\psi}$ and $\hh_{n}^{\phi,\psi}$ have positive sign. And so the polynomial $\hh_{n+1}^{\phi,\psi}$ has also zeros in the intervals $(-\infty,\zeta_{1})$ and $(\zeta_n,+\infty)$. This completes the proof.

\end{proof}
A Turán type inequality is an inequality of the form
$$
p^2(x)-q(x)r(x)>0,\quad x\in I,
$$
where $p,q,r:I\to \RR$ are real functions defined in an interval $I$ of the real line.

It was found by Pál Turán for three consecutive Legendre polynomials ($I=(-1,1)$) and first published by Szegö (\cite{Szet}). It is also true for many other sequences of polynomials, including the ultaspherical, Laguerre and Hermite polynomials (see \cite{Szet}).

As a consequence of Theorem \ref{pep}, we deduce that under mild conditions on the parameters $\phi$ and $\psi$, three consecutive generalized Hermite polynomials $\hh_{j}^{\phi,\psi}$, $j=n-2,n-1,n$, also satisfy a Turán type inequality.

\begin{Theo}\label{ttur}
Let $\phi$ and $\psi$ be two sequences of real numbers with $\psi_i>-2$, $i\ge 1$. Assume that $\phi_{n-1}=\phi_n$ and $\psi_{n-1}=\psi_n$ for some $n$, then
\begin{equation}\label{tur}
(\hh_{n-1}^{\phi,\psi})^2(x)-\hh_{n}^{\phi,\psi}(x)\hh_{n-2}^{\phi,\psi}(x)>0,\quad x\in \RR.
\end{equation}
\end{Theo}

\begin{proof}
Write
\begin{equation}\label{tur1}
r(x)=(\hh_{n-1}^{\phi,\psi})^2(x)-\hh_{n}^{\phi,\psi}(x)\hh_{n-2}^{\phi,\psi}(x)=\left|\begin{matrix} \hh_{n-1}^{\phi,\psi}(x) & \hh_{n}^{\phi,\psi}(x)\\\hh_{n-2}^{\phi,\psi}(x)&\hh_{n-1}^{\phi,\psi}(x)\end{matrix}\right|.
\end{equation}
The polynomial $r$ has then degree at most $2n-2$. Since $\psi_{n-1}=\psi_n$ and $\phi_{n-1}=\phi_n$, it is easy to see that actually $r$ has degree $2n-4$ with leading coefficient equal to
$(\psi_n+2)\prod_{j=1}^{n-2}(\psi_i+2)^2>0$. Hence, if (\ref{tur}) does not hold there will exist $x_0\in \RR$ such that $r(x_0)=0$. And so, there exist $a,b\in \RR$, at least one of them not equal to zero, such that
the polynomials
$$
p(x)=a\hh_{n-1}^{\phi,\psi}(x)+b\hh_{n}^{\phi,\psi}(x),\quad q(x)=a\hh_{n-2}^{\phi,\psi}(x)+b\hh_{n-1}^{\phi,\psi}(x),
$$
have a common zero at $x=x_0$.

On the one hand, since the polynomials $\hh_{n-1}^{\phi,\psi}(x)$ and $\hh_{n-2}^{\phi,\psi}(x)$ interlace their zeros, we have that $q$ has $n-1$ real and simple zeros (see Theorem \ref{obre}).
On the other hand, since $\phi_{n-1}=\phi_n$ and $\psi_{n-1}=\psi_n$, it is easy to see that
$$
p(x)=\Lambda q(x)+(\phi_{n}+x\psi_n)q(x)=-q'(x)+((2+\psi_n)x+\phi_n)q(x).
$$
If $p(x_0)=q(x_0)=0$, then $x_0$ would be a zero of $q$ of multiplicity larger than 1, which it is a contradiction.

\end{proof}

Theorem \ref{ttur} may fail if $\psi_n\not =\psi_{n-1}$, in both cases: when $\phi_n\not =\phi_{n-1}$ or when $\phi_n=\phi_{n-1}$. For instance, for $n=2$ and $\psi_1=1$, the inequality (\ref{tur}) is never true for any real numbers $\psi_2,\phi_1, \phi_2$ as long as $\psi_2>1$.

If $\psi_n =\psi_{n-1}$ and $\phi_n\not =\phi_{n-1}$
then (\ref{tur}) is never true because $r$ (\ref{tur1}) is a polynomial of odd degree $2n-3$ with leading coefficient equal to
$$
(\phi_{n-1}-\phi_n)(\psi_n+2)\prod_{j=1}^{n-2}(\psi_i+2)^2\not =0.
$$

\bigskip

The case $\psi_i=0$ of the generalized Hermite polynomials (\ref{dhg})  is specially interesting. We simplify the notation writing
$$
\hh_n^{\phi,0}(x)=\hh_n^{\phi}(x),\quad n\ge 0.
$$
In order to study it in detail, we need some notation. For $l\ge 1$, we denote $\phi^{\{ l\}}$ for the sequence
\begin{equation}\label{que}
\phi^{\{ l\} }_i=\begin{cases} \phi_i, &1\le i\le l-1,\\
\phi_{i+1}, &l\le i,\end{cases}
\end{equation}
that is, $\phi^{\{ l\} }$ is que sequence obtained by removing the term $\phi_l$ from $\phi$.
First of all, we show that for $\psi_i=0$, the polynomials $\hh_n^{\phi}$, $n\ge 0$,  have the following alternative definition. For $n\ge0$, set
\begin{equation}\label{dpm}
\Phi_i^n\equiv \Phi_i^n(\phi)=\begin{cases} 1, &i=0,\\
\displaystyle \sum_{1\le j_1<\dots<j_i\le n}\phi_{j_1}\cdots\phi_{j_i}, &1\le i\le n,\end{cases}
\end{equation}
so that
\begin{equation}\label{dhpc}
\prod_{i=1}^n(x+\phi_i)=\sum_{j=0}^n\Phi_{n-j}^nx^j.
\end{equation}

\begin{Lem}\label{lrch}
For $n\ge 0$, we have
\begin{equation}\label{dpbg}
\hh_n^{\phi}(x)=\sum_{j=0}^n\Phi_{n-j}^nH_j(x).
\end{equation}
Moreover, for all $l\ge 1$ and $n\ge l-1$, we have
\begin{equation}\label{rrgpl}
\hh_{n+1}^\phi(x)=\Lambda\hh_n^{\phi^{\{ l\}}}(x)+\phi_{l}\hh_n^{\phi^{\{ l\}}}(x).
\end{equation}\end{Lem}

The identity (\ref{dpbg}) shows that the polynomial $\hh_n^{\phi}(x)$ has degree $n$, leading coefficient equal to $2^n$, only depends on the numbers $\phi_1,\dots, \phi_n$ and, moreover, $\hh_n^\phi$ is a symmetric function of $\phi_1,\dots, \phi_n$.

\begin{proof}
If we denote
$$
\Phi_i^{n,l}=\Phi_i^n(\phi^{\{ l\}})
$$
(see (\ref{que})), it is easy to see (from the definition (\ref{dpm})) that for $n\ge l-1$
\begin{equation}\label{xxx}
\Phi_i^{n+1}=\begin{cases} \Phi_0^{n+1}, &i=0,\\
\Phi_i^{n,l}+\phi_l\Phi_{i-1}^{n,l}, &1\le i\le n-1,\\
\phi_l\Phi_{n-1}^{n,l}.\end{cases}
\end{equation}
The identity (\ref{dpbg}) follows now easily by induction on $n$ using (\ref{xxx}) for $l=n+1$.

An easy computation using (\ref{xxx}), (\ref{dpbg}) and (\ref{rbel}) gives
\begin{align*}
\hh_{n+1}^\phi(x)&=\sum_{j=0}^{n+1}\Phi_{n+1-j}^{n+1}H_j(x)\\
&=\phi_l\hh_{n}^{\phi^{\{l\}}}(x)+\sum_{j=0}^{n}\Phi_{n-j}^{n,l}H_{j+1}(x)\\
&=\phi_l\hh_{n}^{\phi^{\{l\}}}(x)+\sum_{j=0}^{n}\Phi_{n-j}^{n,l}\Lambda H_{j}(x)\\
&=\phi_l\hh_{n}^{\phi^{\{l\}}}(x)+\Lambda \left[\sum_{j=0}^{n}\Phi_{n-j}^{n,l}H_{j}(x)\right]\\
&=\phi_l\hh_{n}^{\phi^{\{l\}}}(x)+\Lambda \hh_{n}^{\phi^{\{l\}}}(x).
\end{align*}

\end{proof}

For a positive integer $l$ and a real number $M$ write $\phi^{l,M}$ for the sequence
\begin{equation}\label{sph}
\phi^{l,M}_i=\phi_i+M\delta_{i,l}
\end{equation}
(where $\delta_{i,l}$ denotes de Kronecker delta).

\bigskip

\begin{Theo}\label{pcgb}  Let $\phi$ be a sequence of real numbers.
\begin{enumerate}
\item The polynomial $\hh_n^\phi$ has $n$ real and simple zeros, and for all $l\ge 1$  the zeros of $\hh_{n+1}^\phi$ interlace  the zeros of $\hh_{n}^{\phi^{\{l\}}}$.
\item We denote $\zeta_k(n,\phi)$, $1\le k\le n$, the $k$-th zero of $\hh_n^\phi$, arranging the zeros in increasing order (to simplify the notation and when the context allows it we sometimes will write $\zeta_k, \zeta_k(n)$ or $\zeta_k(\phi)$). Then $\zeta_k(\phi)$ is a decreasing function of $\phi$. More precisely,
we say that  $\phi\le \rho$ if $\phi_j\le \rho_j$ for all $j\ge 1$. Then, if $\phi\le \rho$  we have $\zeta_k(\rho)\le \zeta_k(\phi)$.
\item For a positive integer $l$ and a real number $M\not =0$, consider the sequence $\phi^{l,M}$ (see (\ref{sph})). For $l\le n$, if $M>0$ then the zeros of $\hh_{n}^{\phi^{l,M}}$ interlace the zeros of $\hh_n^\phi$, and if $M<0$ the zeros of $\hh_n^\phi$ interlace the zeros of $\hh_{n}^{\phi^{l,M}}$.
\end{enumerate}
\end{Theo}

\begin{proof}

The proof of the Part (1) is just as that of Theorem \ref{pep}, but using the identity (\ref{rrgpl}) instead of (\ref{dhg}) for $n\ge l-1$.

\medskip

We next prove  the Part (2).

Since $\hh_n^\phi$ only depends on $\phi_i$, $1\le i\le n$, we have that $\zeta_k$ is an smooth function of each $\phi_i$. In order to prove the Part (2), it is enough to prove that $\partial \zeta_k(\phi)/\partial \phi_i<0$, $1\le i\le n$. To simplify the notation, we write $\zeta=\zeta_k$.
Since $\hh_n^\phi(\zeta)=0$, by deriving with respect to $\phi_i$, we deduce
$$
\frac{d \hh_n^\phi(x)}{dx}|_{x=\zeta}\frac{\partial \zeta(\phi)}{\partial \phi_i}+\frac{\partial \hh_n^\phi(x)}{\partial \phi_i}|_{x=\zeta}=0.
$$
Since $\hh_n^\phi$ has simple zeros and $\sign(\lim_{x\to -\infty}\hh_n^\phi(x))=(-1)^n$, it follows that
\begin{equation}\label{sig}
\sign \frac{\partial \zeta(\phi)}{\partial \phi_i}=(-1)^{n+k+1}\sign \left(\frac{\partial \hh_n^\phi(x)}{\partial \phi_i}|_{x=\zeta}\right).
\end{equation}
A simple computation shows that (see (\ref{dpm}))
$$
\frac{\partial \Phi^n_l(\phi)}{\partial \phi_i}=\Phi^{n-1}_{l-1}(\phi^{\{i\}}).
$$
Hence, from (\ref{dpbg}), we deduce
\begin{align*}
\frac{\partial \hh_n^\phi(x)}{\partial \phi_i}&=\sum_{j=0}^n\frac{\partial \Phi_{n-j}^n (\phi)}{\partial \phi_i}\hh_j^\phi(x)\\
&=\sum_{j=0}^{n-1}\Phi_{n-j}^{n-1} (\phi^{\{i\}})\hh_j^\phi(x)=\hh_{n-1}^{\phi^{\{i\}}}(x).
\end{align*}
Using the Part (1) of this Theorem, we have that
$$
\sign \left(\frac{\partial \hh_n^\phi(x)}{\partial \phi_i}|_{x=\zeta}\right)=\sign \hh_{n-1}^{\phi^{\{i\}}}(\zeta)=(-1)^{n+k}.
$$
Hence, using (\ref{sig}), we finally find
$$
\sign \frac{\partial \zeta(\phi)}{\partial \phi_i}=-1.
$$

\medskip

We finally prove the Part (3). A simple computation gives
$$
\Phi_i^n(\phi^{l,M})=\Phi_i^n(\phi)+M\Phi_{i-1}^{n-1}(\phi^{\{l\}}),
$$
for $1\le i\le n$ and $l\le n$, and so
\begin{equation}\label{mcc}
\hh_{n}^{\phi^{l,M}}(x)=\hh_n^\phi(x)+M\hh_{n-1}^{\phi^{\{l\}}}(x).
\end{equation}
Write $\zeta_i$ for the zeros of $\hh_n^\phi(z)$, so that
$$
\zeta_1<\dots<\zeta_n.
$$
Hence $\hh_{n}^{\phi^{l,M}}(\zeta_i)=M\hh_{n-1}^{\phi^{\{l\}}}(\zeta_i)$. It is now enough to take into account the Part (1).

\end{proof}

\medskip

\begin{Rem}\label{mher}
Given a multi-index $\vec{n}=(n_1,\dots, n_r)\in \NN^r$ and the $r$-tuple of real numbers $\vec{c}=(c_1,\dots, c_r)$, $c_i\not =c_j$, $i\not =j$, the multiple Hermite polynomials are defined by the Rodrigues formula (\cite[\S 23.5]{Ism})
\begin{equation}\label{mherd}
H_{\vec{n}}^{\vec{c}}(x)=(-1)^{|\vec{n}|}e^{x^2}\prod_{j=1}^r\left(e^{-c_jx}\frac{d^{n_j}}{dx^{n_j}}e^{c_jx}\right)e^{-x^2},
\end{equation}
where $|\vec{n}|=\sum_{j=1}^r n_j$.
They satisfy the orthogonality conditions
$$
\int_\RR H_{\vec{n}}^{\vec{c}}(x)e^{-x^2+c_jx}x^k=0,\quad k=0,1,\dots, n_j-1,
$$
for $j=1,2,\dots, r$ (so, according to the usual definition of multiple orthogonal polynomials, see \cite{MW}, they are the \textit{type II} multiple Hermite polynomials).

We will prove that the class of all generalized Hermite polynomials $\hh_n^\phi$ is the same class as that of all  multiple Hermite polynomials $H_{\vec{n}}^{\vec{c}}$. More precisely,
\begin{equation}\label{cgblm}
H_{\vec{n}}^{\vec{c}}(x)=\hh_{|\vec{n}|}^{\phi^{{\vec c,\vec{n}}}}(x),
\end{equation}
where
\begin{equation}\label{anph}
\phi^{\vec c,\vec{n}}_i=\begin{cases} -c_1,&i=1,\dots,n_1,\\
-c_2,&i=n_1+1,\dots,n_1+n_2,\\
\dots &\\
-c_r,&i=n_1+\cdots+n_{r-1}+1,\dots,|\vec{n}|.
\end{cases}
\end{equation}
Let us notice that the sequence $\phi^{{\vec c},\vec{n}}$ depends on both the parameters $\vec{c}$ and the multi-index $\vec{n}$.

Multiple Hermite polynomials (\ref{mherd}) satisfy the following two identities. The recurrence relations (\cite[\S 23.5]{Ism})
$$
xH_{\vec{n}}^{\vec{c}}(x)=\frac{1}{2}H_{\vec{n}+\vec{e_k}}^{\vec{c}}(x)+\frac{c_k}{2}H_{\vec{n}}^{\vec{c}}(x)+\sum_{j=1}^r n_jH_{\vec{n}-\vec{e_j}}^{\vec{c}}(x),\quad 1\le k\le r,
$$
where $\vec{e_1}=(1,0,0,\dots, 0)$, $\vec{e_2}=(0,1,0,\dots, 0)$, \dots, $\vec{e_r}=(0,0,0,\dots, 1)$ (notice that we are using a different normalization that in \cite{Ism}).
And the lowering operator (\cite[Eq. (23.8.6)]{Ism}
$$
(H_{\vec{n}}^{\vec{c}})'(x)=2\sum_{j=1}^r n_jH_{\vec{n}-\vec{e_j}}^{\vec{c}}(x).
$$
Combining both identities, we have
$$
H_{\vec{n}+\vec{e_k}}^{\vec{c}}(x)=2xH_{\vec{n}}^{\vec{c}}(x)-(H_{\vec{n}}^{\vec{c}})'(x)-c_kH_{\vec{n}}^{\vec{c}}(x).
$$
The identity (\ref{cgblm}) can then be proved easily from the definition (\ref{dhg}) ($\psi=0$) using induction on $|\vec{n}|$.

In particular, Part (2) of Theorem \ref{pcgb} implies that the $k$-th zero of the multiple Hermite polynomial $H_{\vec{n}}^{\vec{c}}$, $1\le k\le |\vec{n}|$, is an increasing function of the parameters $c_1,\dots, c_r$.
\end{Rem}

\section{Zeros of linear combinations of finitely many Hermite polynomials}\label{secH2}

If we fix a nonnegative integer $K$ and real numbers $\gamma_j$, $j=0,\dots, K$, with $\gamma_0=1$, $\gamma_K\not =0$, in this section we return to the problem of determining the number of the real zeros of the polynomials
\begin{equation}\label{polqn}
q_n(x)=\sum_{j=0}^K\gamma_jH_{n-j}(x),\quad n\ge K.
\end{equation}
Since $\Lambda H_n=H_{n+1}$, we have that for $n\ge K$, $q_n=\Lambda^{n-K}q_K$.

We associate to the real numbers $\gamma_j$, $0\le j\le K$, the polynomial
\begin{equation}\label{polp}
P(x)=\sum_{j=0}^K\gamma_jx^{K-j}.
\end{equation}

On the one hand, the orthogonality of the Hermite polynomials implies that the polynomials $q_n$, $n\ge K$, has at least $n-K$ zeros (see \cite[Lemma 3.1]{Dur0}), and on the other hand, Corollary \ref{1her} implies that
for $n$ big enough (depending on the $\gamma_j$'s, and hence on $P$) all the zeros of $q_n$ are real. However, Theorem \ref{pcgb} says quite more about the real zeros of $q_n$: they are always real and simple for $n\ge K$ if all the zeros of $P$ are real, and when some of the zeros of $P$ are not real then the \textit{big enough} (mentioned above) only depends of the non real zeros of the polynomial $P$.
We next prove the following Corollary (which it is a more detailed version of Theorem \ref{pcgb}).

\begin{Cor}\label{2her} Let $\theta_i$, $i=1,\dots, 2m$, the non real zeros of the polynomial $P$ (\ref{polp}), and define
$$
P_{\nr}(x)=\prod_{j=1}^{2m}(x-\theta_j)=\sum_{j=0}^{2m}\tau_jx^{2m-j},
$$
where $\tau_j$ are real numbers. Let $n_0$ be the first positive integer such that the polynomial
$$
\sum_{j=0}^{2m}\tau_jH_{n_0-j}(x)
$$
has only real and simple zeros (Corollary \ref{1her} implies that such positive integer always exists and
$n_0\le \max\{(2m-1)^24^{2m-2}\max^2\{|\tau_j|,2\le j\le 2m\},4m\}$). Then the polynomial $q_n$ (\ref{polqn})
has only real zeros for $n\ge n_0+K-m$, and the zeros of $q_{n+1}$ interlace the zeros of $q_n$. Moreover, if $P$ has only real zeros then $n_0=0$.
\end{Cor}

\begin{proof}
We start proving the case when $P$ has only real zeros (i.e., $m=0$). Write $\theta_i$, $i=1,\dots ,K$, for them. We have using Lemma \ref{lrch} that $q_n=\hh_n^\phi$, $n\ge K$, where $\phi_i=-\theta_i$, $1\le i\le K$, and $\phi_i=0$, $i\ge K+1$. And so Theorem \ref{pcgb} says that $q_n$ has only real zeros for $n\ge K$, and the zeros of $q_{n+1}$ interlace the zeros of $q_n$.

We next prove the case when $m>0$.

The proof of the case $K=2m$ is as follows. Since $K=2m$ we have $P_{\nr}=P$ and then  $\gamma_j=\tau_j$, $0\le j\le 2m$.
The hypothesis and Corollary \ref{1her} say that the polynomial $q_{n_0}$ has only real and simple zeros. Since $\Lambda H_n=H_{n+1}$, we have that $q_n=\Lambda^{n-n_0}(q_{n_0})$ for $n\ge n_0$. It is then enough to use Lemma \ref{hoiz} (which it also proves that the zeros of $q_{n+1}$  interlace the zeros of $q_n$).

If $K-2m=1$, let $\theta$ be the real zero of $P$. Denoting $B_j=\tau_j$, $A_j=\gamma_j$ and using the notation of Lemma \ref{lem2},
we have (see (\ref{idco}))
$$
q_{n+1}(x)=q_{n+1}^A(x)=q_{n+1}^{B}(x)-\theta q_{n}^{B}(x),
$$
and since $q_n^{B}=q_n^{\tau}$, and we have already proved that the zeros of $q_{n+1}^{\tau}$ interlace the zeros of $q_n^{\tau}$, $n\ge n_0$, Lemma \ref{enze} gives that $q_n$ has only real zeros for $n\ge n_0+1$.
Since $\Lambda q_n=q_{n+1}$, Lemma \ref{hoiz} implies that the zeros of $q_{n+1}$  interlace the zeros of $q_n$ for $n\ge n_0+1$.

The cases $K-2m\ge 2$ can be proved similarly.

\end{proof}

As a consequence of Corollary \ref{2her}, we deduce the following Turán type inequality. We omit the proof because is essentially the same as that of Theorem \ref{ttur} (using that $\Lambda q_n=q_{n+1}$).

\begin{Cor}\label{3her} The polynomial $q_n$ (\ref{polqn}) satisfy the following Turán type inequality:
$$
q_{n+1}^2(x)-q_{n+2}(x)q_{n}(x)>0,\quad x\in \RR,
$$
where $n\ge K$ if all the zeros of $P$ are real, and $n$ has to be taken big enough if $P$ has non real zeros (the big enough depends only on the non real zeros of $P$ as explain in Corollary \ref{2her}).
\end{Cor}

\section{Other normalization of the Hermite polynomials}\label{secH3}
We next take the following normalization of the Hermite polynomials
\begin{equation}\label{hen}
\Hh_n(x)=\frac{1}{2^nn!}H_n(x),
\end{equation}
and consider linear combinations of the polynomials $\Hh_n$ of the form
\begin{equation}\label{polqns}
q_n(x)=\sum_{j=0}^K\gamma_j\Hh_{n-j}(x),\quad n\ge 0,
\end{equation}
where $K$ is a positive integer and
$\gamma_j$, $j=0,\dots, K$, are real numbers with $\gamma_0=1$, $\gamma_K\not =0$ (we take $\Hh_j=0$, for $j<0$).

Since $\Hh_n(x)=\hat H_n(x)/n!$ (where $\hat H_n$ denotes the monic Hermite polynomials), we can not use Theorem 3.3 in \cite{Dur0}
to study the zeros of the polynomials $q_n$ (because the bounds (3.9) in that Theorem never hold for any sequence $(\tau_n)_n$ satisfying $\lim_n\tau_n=0$).

However, using another approach we can completely describe the zeros of the polynomials $q_n$. In this case, the key is that the polynomials $q_n$ are Appell polynomials:
$q_n'=q_{n-1}$, $n\ge 1$. Moreover, the generating function  for the Hermite polynomials (see \cite[Identity (9.15.10)]{KLS}) gives
\begin{equation}\label{gen}
\sum_{n=0}^\infty q_n(x)z^n=e^{xz-z^2/4}R(z),
\end{equation}
where
\begin{equation}\label{psov}
R(x)=\sum_{j=0}^K\gamma_jx^{j}.
\end{equation}
Write
\begin{equation}\label{psov1}
P(x)=\sum_{j=0}^K\gamma_jx^{K-j},
\end{equation}
so that $R(x)=x^KP(1/x)$. As a consequence of (\ref{gen}) and (\ref{psov1}), the polynomials $(q_n)_n$ enjoy the following asymptotic (see \cite[Theorem 1.1]{drh})
\begin{equation}\label{asy}
\lim_n\left(\frac{z}{n+1}\right)^nn!q_n\left(\frac{n+1}z\right)=z^Ke^{-z^2/4}P(1/z),
\end{equation}
uniformly in compact sets of the complex plane.

We are now ready to prove the Theorem \ref{zeri} in the Introduction which describes the structure of the zeros of the polynomials $q_n$ (\ref{polqns}) (let us note that the polynomials $P$ (\ref{psov1}) and $R$ (\ref{psov}) have the same number of positive, negative and real zeros, respectively).

\begin{proof}[Proof of the Theorem \ref{zeri}]
Write $Z^{\nr}(n)$ for the number of non real zeros of $q_n$. Since $q_n=q'_{n+1}$, we deduce that
\begin{equation}\label{znd}
\mbox{$Z^{\nr}(n)$ is a non-decreasing function of $n$}.
\end{equation}

We proceed by induction on $K-N^{\nr}$.

If $N^{\nr}=K$, the asymptotic (\ref{asy}) implies that there exists $n_0$ (which we take it to be the smallest one) such that $q_{n_0}$ has at least $K$ non real zeros. Using (\ref{znd}), we deduce that $q_n$ has also at least $K$ non real zeros for $n\ge n_0$. \cite[Lemma 3.1]{Dur0} then implies that $q_n$ has exactly $n-K$ real zeros and they are simple. This also says that $n_0$ is the smallest positive integer such that $q_{n_0}$ has exactly $K$ non real zeros. Since $q_{n+1}'=q_n$, the real zeros of $q_{n+1}$ interlace the real zeros of $q_n$. This proves Part (2) of the Theorem \ref{zeri} for $N^{\nr}=K$.

Assume next that $K-N^{\nr}>1$. Then we have that $P$ has at least a real zero $\theta$.
Write, as in the Lemma \ref{lem2},
\begin{equation}\label{poq}
P_{B}(x)=P(x)/(x-\theta)=\sum_{j=0}^{K-1}B_jx^{K-1-j},
\end{equation}
and define $q_n^B(x)=\sum_{j=0}^{K-1}B_j\Hh_{n-j}(x)$ (\ref{pidco}) (with the notation in (\ref{pidco}) $q_n=q_n^{A;\theta}$).
The induction hypothesis says that there exists $n_1$, the smallest positive integer such that $q_{n_1}^B$ has exactly $N^{\nr}$ non real zeros, and that for $n\ge n_1$ the polynomials $q_{n+1}^B$ and $q_{n}^B$ have exactly $n-N^{\nr}+1$ and $n-N^{\nr}$ real and simple zeros, respectively, and the real zeros of $q_{n+1}^B$ interlace the real zeros of $q_{n}^B$.
Hence, if we write $\zeta_1<\dots <\zeta_{n-N^{\nr}}$ for the real zeros of $q_{n}^B$, we deduce that $q_{n+1}^B$ has exactly one zero in each interval
$(\zeta_i,\zeta_{i+1})$, $0\le i\le n-N^{\nr}$, where $\zeta_0=-\infty$ and $\zeta_{n-N^{\nr}+1}=+\infty$.

Using (\ref{idco}) we have
\begin{equation}\label{ftqn}
q_{n+1}(x)=q_{n+1}^B(x)-\theta q_{n}^B(x).
\end{equation}
This gives
$$
q_{n+1}(\zeta_i)q_{n+1}(\zeta_{i+1})=q_{n+1}^B(\zeta_i)q_{n+1}^B(\zeta_{i+1}),
$$
and we conclude that $q_{n+1}$ has also at least $n-N^{\nr}+1$ real zeros for $n\ge n_1$. In particular $Z^{\nr}(n_1+1)\le N^{\nr}$.
The asymptotic (\ref{asy}) implies that there exists $n_0$ (which we take it to be the smallest one) such that $q_{n_0}$ has at least $N^{\nr}$ non real zeros, and hence $N^{\nr}\le Z^{\nr}(n_0)$. We then deduce that
$n_0\ge n_1+1$ (because of (\ref{znd})). We also have that then $q_{n_0}$ has exactly $N^{\nr}$ non real zeros and $q_{n}$ has exactly $n-N^{\nr}$ real and simple zeros for $n\ge n_0$.
Moreover, we have also proved that the zeros of $q_{n+1}$ interlace the zeros of $q_{n}^B$. And they also interlace the zeros of $q_{n}$, again because $q_{n}=q_{n+1}'$.

If $N^{\nr}=0$, then all the zeros of $q_n$ has to be real for $n\ge 0$, because $q_n=q_{n+1}'$. This completes the proof of the Theorem.

\end{proof}

\medskip

Part (2) of the Theorem \ref{zeri} can be completed as follows.

\begin{Cor}\label{2her2} Assume in Theorem \ref{zeri} that $N^{\nr}>0$. For a real number $\theta$, define
$$
q_{n,\theta}(x)=q_{n}(x)-\theta q_{n-1}(x).
$$
Then there exists a nonnegative integer $n_*$, which does not depend on $\theta$,
such that for $n\ge n_*$, the polynomial $q_{n,\theta}$  has exactly $n-N^{\nr}$ real and simple zeros and $N^{\nr}$ non real zeros, and the real zeros of $q_{n+1,\theta}$ interlace the real zeros of $q_{n,\theta}$.
\end{Cor}

\begin{proof}
Let $n_0$ be as in Part (2) of Theorem \ref{zeri}, i.e., the smallest positive integer such that $q_{n_0}$ has $N^{\nr}$ non real zeros.

Following the notation of Lemma \ref{lem2}, define
\begin{align}\label{poq2}
P_A(x)&=(x-\theta)P(x)=\sum_{j=0}^{K+1}A_jx^{K+1-j},\\\label{poq3}
q_n^A(x)&=\sum_{j=0}^{K+1}A_j\Hh_{n-j}(x),
\end{align}
and write $P_B(x)=P(x)$, and $q_n^B(x)=q_n(x)$, so that (see (\ref{idco}))
\begin{equation}\label{ftqn2}
q_{n}^A(x)=q_{n}^B(x)-\theta q_{n-1}^B(x)=q_{n,\theta}(x).
\end{equation}
Define finally
\begin{align*}
f_n(z)&=\left(\frac{z}{n+1}\right)^nn!q_n^B\left(\frac{n+1}z\right),\\
g_n(z)&=\left(\frac{z}{n+1}\right)^nn!q_{n-1}^B\left(\frac{n+1}z\right),\\
f(z)&=z^{K}e^{-z^2/4}P(1/z).
\end{align*}
The asymptotic (\ref{asy}) for $q_n^B$ gives
\begin{equation}\label{asy2}
\lim_nf_n(z)=f(z).
\end{equation}
In turns, from  the asymptotics (\ref{asy}) for $q_n^A$ and (\ref{asy2}) for $q_n^B$ we deduce
$$
\lim_ng_n(z)=zf(z).
$$
Lemma \ref{mamon} guarantees the existence of a positive integer $n_*$ which does not depend on $\theta$, and which can be taken $n_*\ge n_0$, such that for $n\ge n_*$
$$
f_n(z)-\theta g_n(z)=\left(\frac{z}{n+1}\right)^nn!q_n^A\left(\frac{n+1}z\right)
$$
has at least $N^{\nr}$ non real zeros. Hence $q_n^A$ has at least $N^{\nr}$ non real zeros for $n\ge n_*$. Since $n_*\ge n_0$, Theorem \ref{zeri} says that $q^B_n$ has exactly $n-N^{\nr}$ real and simple zeros for $n\ge n_*$, and the real zeros of $q_{n+1}^B$ interlace the zeros of $q_{n}^B$. Proceeding as in the proof of Theorem \ref{zeri} (using (\ref{ftqn2})), we can deduce that $q^A_n$ has exactly $n-N^{\nr}$ real and simple zeros for $n\ge n_*$, and the real zeros of $q_{n+1}^A$ interlace the zeros of $q_{n}^A$.
\end{proof}

\bigskip

As a consequence of Theorem \ref{zeri} and Corollary \ref{2her2}, we deduce the following Turán type inequality.

\begin{Cor} The polynomial $q_n$ (\ref{polqns}) satisfy the following Turán type inequality:
\begin{equation}\label{tur2}
q_{n-1}^2(x)-q_n(x)q_{n-2}(x)>0,\quad x\in \RR,
\end{equation}
where $n\ge 2$ if all the zeros of $P$ are real, and $n$ has to be taken big enough if $P$ has non real zeros.
\end{Cor}

\begin{proof}
Write, as in the proof of Theorem \ref{ttur},
\begin{equation}\label{tur3}
r(x)=q_{n-1}^2(x)-q_{n}(x)q_{n-2}(x)=\left|\begin{matrix} q_{n-1}(x) &q_{n}(x)\\ q_{n-2}(x)&q_{n-1}(x)\end{matrix}\right|.
\end{equation}
A simple computation shows that polynomial $r$ has then degree $2n-2$ with leading coefficient equal to
$1/((n-1)!n!)$.

Take $n_*=2$, if $P$ has only real zeros, and $n_*$ as in Corollary \ref{2her2} if $P$ has some non real zeros.

Hence, if (\ref{tur2}) does not hold for $n\ge n_*$, there will exist $x_0\in \RR$ such that $r(x_0)=0$. And so, there exist $a,b\in \RR$, at least one of them not equal to zero, such that
the polynomials
$$
p(x)=aq_{n-1}(x)+bq_{n}(x),\quad q(x)=aq_{n-2}(x)+bq_{n-1}(x),
$$
have a common zero at $x=x_0$. Since $p'=q$, that means that $p$ has a zero at $x=x_0$ of multiplicity at least $2$.
Since the real zeros of $q_n$ are simple for $n\ge n_*$, we can assume that $b\not =0$.
This implies that the polynomial $q_n-\theta q_{n-1}$ has a zero at $x=x_0$ of multiplicity at least $2$, where $\theta=-a/b$. This contradicts either Theorem \ref{zeri} (if $P$ has only real zeros) or
Corollary \ref{2her2} (if $P$ has some non real zeros).
\end{proof}

\bigskip

We next display the asymptotic behaviour of the zeros of the polynomials $q_n$ (\ref{polqns}).

\begin{Cor}\label{cas2} Assume that the polynomial $P$ (\ref{psov1}) has $N^{\nr}$ non real zeros and $N^-$ negative zeros. For $n$ big enough write $\zeta_j(n)$, $1\le j\le n-N^{\nr}$, for the real zeros of the polynomial $q_n$ (\ref{polqn}) arranged in increasing order, and $\zeta_j^{\nr}(n)$, $1\le j\le N^{\nr}$, for the non  real zeros of the polynomial $q_n$ arranged according to the complex lexicographic order. Similarly, write $\theta_j$, $1\le j\le K-N^{\nr}$, for the real zeros of $P$ (arranged also in increasing order) and $\theta_j^{\nr}(n)$, $1\le j\le N^{\nr}$, for the non  real zeros of the polynomial $P$ (arranged according to the complex lexicographic order).
\begin{enumerate}
\item Asymptotic for the central real zeros: for $j\in\ZZ$,
\begin{align}\label{lec9}
\lim_n \sqrt {2n}\zeta_{j+[(n-K)/2]+N^-+1}(n)&=\begin{cases}\frac{\pi}2+j\pi,& \hbox{$n-K$ is even}\\
j\pi,&\hbox{$n-K$ is odd}.\end{cases}
\end{align}
\item Asymptotic for the leftmost and rightmost real zeros:
\begin{align}\label{lec7}
\lim_n \frac{\zeta_{j}(n)}{n+1}&=\theta_j,\quad 1\le j\le N^-,\\\label{lec8}
\lim_n \frac{\zeta_{n-K+j}(n)}{n+1}&=\theta_{j},\quad N^-+1\le j\le K-N^{\nr}.
\end{align}
\item Asymptotic for the non real zeros:
$$
\lim_n \frac{\zeta^{\nr}_{j}(n)}{n+1}=\theta_{j}^{\nr},\quad 1\le j\le N^{\nr}.
$$
\item For a bounded continuous function $f$ in $\RR$, we have
$$
\lim_n\frac{1}{n}\sum_{j=1}^{n-N^{\nr}}f\left(\frac{\zeta_j(n)}{\sqrt{2n}}\right)=\frac{2}{\pi}\int_{-1}^1 f(x)\sqrt{1-x^2}dx.
$$
\end{enumerate}
\end{Cor}

\begin{proof}
Consider next the well-known Mehler-Heine formula for the Hermite polynomials:
\begin{align}\label{lec}
\lim_n\frac{(-1)^n\sqrt{n\pi}}{2^{2n}n!}H_{2n}\left(\frac{x}{2\sqrt n}\right)&=\cos x,\\\label{lec0}
\lim_n\frac{(-1)^n\sqrt{\pi}}{2^{2n+1}n!}H_{2n+1}\left(\frac{x}{2\sqrt n}\right)&=\sin x
\end{align}
(see \cite[Identities 18.11.7 and 18.11.8]{nist}).

This Mehler-Heine formula  easily gives
\begin{align}\label{lec11}
\lim_n\frac{(n-K)!\sqrt{\pi(n-K)/2}}{(-1)^{(n-K)/2}((n-K)/2)!}q_{n}\left(\frac{x}{\sqrt{2(n-K)}}\right)&=\gamma_K\cos x,\quad \hbox{$n-K$ is even}\\\label{lec12}
\lim_n\frac{(-1)^{(n-K-1)/2}(n-K)!\sqrt{\pi}}{((n-K-1)/2)!}q_{n}\left(\frac{x}{\sqrt{2(n-K)}}\right)&=\gamma_k\sin x,\quad \hbox{$n-K$ is odd}
\end{align}

In order to prove the corollary, we use the following Theorem due to Beardon and Driver.

\begin{Theo}\cite[Theorem 3.1]{BDr}\label{BeDr} Let $(p_n)_n$ be orthogonal polynomials with respect to a positive measure. Fix $0<r<n$ and let $\xi_i(n)$, $i=1,\dots, n$, the zeros of $p_n$ listed in increasing order. Let $P$ be a polynomial in the span of $p_r,\dots, p_n$. Then at least $r$ of the intervals $(\xi_i,\xi_{i+1})$ contain a zero of $P$.
\end{Theo}

For each $n$, consider the set of nonnegative integers
$$
I_n=\{j\in \NN: \mbox{$\zeta_j(n)$ is between the zeros of $H_n$}\}.
$$
Theorem \ref{BeDr} says that $I_n$ has at least $n-K$ elements.

The asymptotic (\ref{asy}) implies that for each $i=1, \dots , K-N^{\nr}$ and $n\ge 0$, there exists $1\le j(i,n)\le n-N^{\nr}$ such that
$$
\lim_{n\to \infty}\frac{\zeta_{j(i,n),n}}{n+1}=\theta_i.
$$
If we write $\varsigma_j(n)$, $1\le j\le n$, for the zeros of the Hermite polynomial $H_n$, using the well-known bound
$$
|\varsigma_{j}(n)|\le \sqrt{2n+3}
$$
(see \cite[p.130]{Sze}), we have for $j\in I_n$
$$
\left\vert\frac{\zeta_j(n)}{n+1}\right\vert \le \frac{\sqrt{2n+3}}{n+1}.
$$
We can then conclude that for $n$ big enough $\zeta_{j(i,n),n}\not \in I_n$, $1\le i \le K-N^{\nr}$ (note that $P(0)=\gamma_K\not =0$ and so $\theta_i\not =0$, $1\le i\le K$).

This shows that we can take $j(i,n)=i$, $1\le i \le N^{-}$, and $j(i,n)=n-K+i$, $N^-+1\le i \le K-N^{\nr}$. This proves Part (2) of the corollary.

Part (3) follows as a consequence of the asimptotic (\ref{asy}).

Part (1) is now a consequence of the Mehler-Heine type formula  (\ref{lec11}) and the Hurwitz's Theorem.

The part (4) also follows from Theorem \ref{BeDr}, taking into account the weak scaling limit (\cite{Dei})
\begin{equation}\label{wsl}
\lim_n\frac{1}{n}\sum_{j=1}^nf\left(\frac{\varsigma_j(n)}{\sqrt{2n}}\right)=\frac{2}{\pi}\int_{-1}^1 f(x)\sqrt{1-x^2}dx,
\end{equation}
for the counting measure of the zeros $\varsigma_j(n)$, $1\le j\le n$, of the Hermite polynomials.

\end{proof}

\bigskip

We finish the paper studying the number of positive and negative zeros of the polynomials $q_n$ (\ref{polqns}).

When all the zeros of $P$ (\ref{psov1}) are real we have the following result.

\begin{Cor}\label{ultmm} Assume that all the zeros of the polynomial $P$ (\ref{psov1}) are real, of which $N^-$ are negative. For $n$ big enough, if $n-K$ is even, $q_n$ has $(n-K)/2+N^-$ negative zeros and $(n-K)/2+K-N^{\nr}-N^-$ positive zeros, and if $n-K$ is odd, $q_n$ has at least $(n-K-1)/2+N^-$ negative zeros and at least $(n-K-1)/2+K-N^{\nr}-N^-$ positive zeros.
\end{Cor}

\begin{proof}
Assume first that $n-K$ is even. We next prove that $q_n$ has $(n-K)/2+N^-$ negative zeros and $(n-K)/2+K-N^-$ positive zeros. We proceed by induction on $K$. The case $K=0$ is the Hermite case. Assume $K>1$. Hence $P$ has at least one zero $\theta$.

Write, as in the proof of Theorem \ref{zeri},
$$
P_{B}(x)=P(x)/(x-\theta)=\sum_{j=0}^{K-1}B_jx^{K-1-j},
$$
and define $q_n^B(x)=\sum_{j=0}^{K-1}B_j\Hh_{n-j}(x)$. Hence, we get from (\ref{ftqn})
\begin{equation}\label{ftqnx}
q_{n}(x)=q_{n}^B(x)-\theta q_{n-1}^B(x).
\end{equation}
If  $P$ has a positive zero $\theta>0$, (\ref{ftqnx}) shows that the zeros of $q_n$ interlace the zeros of $q_{n-1}^B$. Since $n-1-(K-1)=n-K$ is even, the induction hypothesis implies that
$q_{n-1}^B$ has $(n-K)/2+N^-$ negative zeros and $(n-K)/2+K- N^--1$ positive zeros. So, $q_n$ has at least $(n-K)/2+N^-$ negative zeros and at least $(n-K)/2+K-N^--1$ positive zeros. We next prove that $q_n$ has one more positive zero. Indeed, write $v=(n-K)/2+N^-$. Then, on the one hand, using (\ref{ftqn}) we have
\begin{equation}\label{ptt1}
\sign(q_n(\zeta_v))=\sign (q_n^B(\zeta_v))=(-1)^{(n-K)/2+v+n}
\end{equation}
(because $v=(n-K)/2+N^-$ and $\zeta_v$ is the $v$-th zero of $q_n^B$).

On the other hand
$$
q_n(0)=\sum_{j=0}^K\gamma_j\Hh_{n-j}(0)=\sum_{\substack{{j=0}\\{\footnotesize \mbox{$n-j$ even}}}}^K\frac{(-1)^{(n-j)/2}}{2^{n-j}((n-j)/2)!}\gamma_j.
$$
Hence, for $n$ big enough we deduce
\begin{equation}\label{ptt2}
\sign q_n(0)=\sign(\gamma_K (-1)^{(n-K)/2})=(-1)^{(n-K)/2+K+N^-}
\end{equation}
(because $P(0)=\gamma_K$ and $P$ has degree $K$ and $N^-$ negative zeros).

Since $n-K$ is even, (\ref{ptt1}) and (\ref{ptt2}) imply that $q_n$ does not vanish in $(\zeta_v,0)$, and hence $q_n$ has to vanish in $(0,\zeta_{v+1})$. This proves that $q_n$ has  $(n-K)/2+N^-$ negative zeros and  $(n-K)/2+K-N^-$ positive zeros.

If $\theta<0$, the proof is similar.

If $n-K$ is odd, then $n-K+1$ is even, since $q_{n-K+1}'=q_{n-K}$, we deduce that $q_n$ has at least $(n-K-1)/2+N^-$ negative zeros and at least $(n-K-1)/2+K-N^-$ positive zeros.
\end{proof}

\medskip

When all the zeros of $P$ are non real, we have the following conjecture:

\noindent \textit{Conjecture}. Assume that all the zeros of the polynomial $P$ (\ref{psov1}) are non real.
For $n$ big enough, if $n-K$ is even, $q_n$ has $(n-K)/2$ negative zeros and $(n-K)/2$ positive zeros (i.e., equal number of positive and negative zeros), and if $n-K$ is odd, $q_n$ has at least $(n-K-1)/2$ negative zeros and at least $(n-K-1)/2$ positive zeros.

If the conjecture in true, then proceeding as in the proof of Corollary \ref{ultmm}, it would follow the following: Assume that the polynomial $P$ (\ref{psov1}) has $N^{\nr}$ non real zeros and $N^-$ negative zeros. For $n$ big enough, if $n-K$ is even, $q_n$ has $(n-K)/2+N^-$ negative zeros and $(n-K)/2+K-N^{\nr}-N^-$ positive zeros, and if $n-K$ is odd, $q_n$ has at least $(n-K-1)/2+N^-$ negative zeros and at least $(n-K-1)/2+K-N^{\nr}-N^-$ positive zeros.

\medskip

%%%%%%%%%%%%%%%%%%%%%%%%%%%%%%%%%%%%%%%%%%%%%%%%%%%%%%%%%%%%%%%%%%%%%%%

%%%%%%%%%%%%%%%%%%%%%%%%%%%%%%%%%%%%%%%%%%%%%%%%%%%%%%%%%%%%%%%%%%%%%%%

\end{document}